\magnification\magstephalf
\baselineskip14pt
\def\pfbox
  {\hbox{\hskip 3pt\lower2pt\vbox{\hrule
  \hbox to 5pt{\vrule height 7pt\hfill\vrule}
  \hrule}}\hskip3pt}
\def\cd#1{\rlap{#1}\raise.2ex\hbox to.5em{\tensy\hss\char"D\hss}}
\def\bib[#1] {\par\noindent\hangindent 20pt\hbox to20pt{[#1]\hfil}}
\def\AW{Addison\kern.1em--Wesley}

\def\meno{\medskip\noindent}
\def\cm{3}
\def\ss{8}
\def\gs{2}
\def\gi{4}
\def\fp{1}
\def\mw{7}
\def\ks{5}
\def\ms{6}

\centerline{\bf An exact analysis of stable allocation}
\medskip
\centerline{Donald E. Knuth}
\centerline{Computer Science Department, Stanford University}

\bigskip
{\narrower\smallskip\noindent
{\bf Abstract.} Shapley and Scarf [\ss] introduced a notion of stable
allocation between traders and indivisible goods, when each trader has
rank-ordered each of the goods. The purpose of this note is to prove that the
distribution of ranks after allocation is the same as the distribution of
search distances in uniform hashing, when the rank-orderings are independent
and uniformly random. Therefore the average sum of final ranks is just
$(n+1)H_n-n$, and the standard deviation is $O(n)$.
 The proof involves a family of  interesting one-to-one correspondences
between permutations of a special kind.
\smallskip}

\meno
{\bf 0. Introduction.}
Suppose $n$ traders have $n$ indivisible goods to trade, and each trader~$k$
has ranked the goods of all traders (including himself) as a permutation
$$p_k=p_{k1} \ p_{k2} \ \ldots \ p_{kn}$$
of $\{1,2,\ldots,n\}$. If $i$ precedes~$j$ in this list, we write ``$i>j \
(k)$'' and say that $k$ prefers~$i$ to~$j$. An allocation of goods to traders
is a permutation $g_1\ldots g_n$ of $\{1,\ldots,n\}$ such that trader~$k$
gets~$g_k$. Shapley and Scarf [\ss] defined what they called a ``core
allocation''~$g$, which is {\it stable\/} in the following sense: If $C$ is any
coalition of traders (a~nonempty subset), and if $h$ is any allocation of the
goods of~$C$ to the members of~$C$, then 
$$h_k\geq g_k\ (k)\quad \hbox{for all $k\in C$}
\quad  \Longrightarrow \quad h_k=g_k\quad\hbox{for all $k\in C$.}\eqno(\ast)$$

For example, suppose $n=3$ and the preference rankings are
$$\eqalign{p_1&=2\,1\,3\,,\cr
p_2&=3\,1\,2\,,\cr
p_3&=2\,3\,1\,.\cr}$$
Then just three of the six possible allocations satisfy $(\ast)$ when $C$ is
the full set $\{1,2,3\}$, namely $2\,1\,3$, $2\,3\,1$, and $1\,3\,2$. But
$2\,1\,3$ is unstable because traders~2 and~3 can both improve their lot by
swapping goods between themselves. Similarly, $2\,3\,1$ is unstable, but it is
not quite as bad: If the coalition $\{2,3\}$ exchanges goods, trader~3 is
happier than he was before, while 2 is no worse off. The remaining allocation,
$1\,3\,2$, is stable.

It is not immediately obvious that a stable allocation always exists, for all
$n!^n$ possible ranking sequences. 
But Shapley and Scarf presented an algorithm by David Gale that always finds
one. In fact, there is always {\it exactly one\/}  stable allocation. Gale's
procedure is similar to the famous Gale-Shapley algorithm for stable marriage
[\gs,\gi], but it incorporates a new twist. 

Alan Frieze and Boris Pittel [\fp] recently analyzed Gale's algorithm and
discovered some remarkable simplifications in the course of their study. One of
the main purposes of the present note is to exhibit some underlying
combinatorial structure that accounts for the surprising phenomena they
discovered. Frieze and Pittel proved, among other things, that the total sum of
ranks in the stable allocation~$g$ (i.e., the sum $r_1+\cdots +r_n$, where
$g_k=p_{kr_k}$) is between $({1\over 2}-\epsilon)n\ln n$ and $(1+\epsilon)n\ln
n$ with high probability as $n\rightarrow\infty$, assuming that the
preferences~$p_k$ are independent and uniformly random.
We will deduce the exact distribution of $r_1+\cdots +r_n$, showing in
particular that its mean value is $(n+1)H_n-n$, and in fact we will see
that the joint distribution of the multiset $\{r_1,\ldots,r_n\}$ has a
particularly simple form.

\meno
{\bf 1. Uniform hashing.}
First let's consider a simpler problem, namely to find an allocation~$g$ that
satisfies~$(\ast)$ when $C$ is the full set $\{1,\ldots,n\}$ but not
necessarily for any other~$C$. Such an allocation is {\it locally optimal}, in
the sense that no trader can improve his selection unless some other trader
loses ground. It's easy to achieve such an allocation by simply letting $g_k$
be the first item of list~$p_k$ that is not in $\{g_1,\ldots,g_{k-1}\}$, for
$k=1,\ldots,n$. 

This trivial allocation algorithm, ``first-come first-served,'' is precisely
the method of {\it uniform hashing\/} that arises in the study of information
retrieval~[\ks], when the preference lists (called ``hash sequences'' in that
context) are randomly chosen. The analysis of uniform hashing is particularly
simple, and we will see below that stable allocation can be reduced to the same
analysis. 

A more general way to obtain a locally optimal allocation is to let $\pi$ be
any permutation of $\{1,\ldots,n\}$, and then to let $g_{\pi(k)}$ be the first
item of $p_{\pi(k)}$ that is not in $\{g_{\pi(1)}\ldots g_{\pi(k-1)}\}$, for
$k=1,\ldots,n$. The permutation~$\pi$ gives top priority to trader $\pi(1)$,
then to $\pi(2)$, and so on. Indeed, {\it every\/} allocation that
satisfies~$(\ast)$ for $C=\{1,\ldots,n\}$ will be found by this method, for
some~$\pi$. The reason is that we must have $g_k=p_{k1}$ for some~$k$, in any
locally optimum~$g$, because the mapping from~$k$ to~$p_{k1}$ for all~$k$
always contains a cycle; any such cycle can be used to improve an allocation in
which no trader has his first choice. Let $\pi(1)$ be any value of~$k$ with
$g_k=p_{k1}$. Remove $k$ from all preference lists and apply the same reasoning
recursively to the remaining $n-1$ traders. This defines a permutation
$\pi(1)\,\pi(2)\,\ldots\,\pi(n)$ such that $g_{\pi(k)}$ is the favorite of
trader~$\pi(k)$ in $\{1,\ldots,n\}\setminus\{g_{\pi(1)},\ldots,g_{\pi(k-1)}\}$.

Since every locally optimum allocation is obtained in this way using
some~$\pi$, the stable allocation must itself result from some~$\pi$. And when
the preference lists are random, any~$\pi$ behaves like any other. Thus we
might expect that stable allocation statistics are essentially identical to the
statistics of uniform hashing. This, in fact, is true, but we must be careful
to make the argument rigorous.

\meno
{\bf 2. An algorithm.}
Let's now consider a simple algorithm that computes the stable
 allocation, given
any sequence of preference rankings $p_1\ldots p_n$. The following procedure is
a sequential variant of Gale's parallel method, analogous to the McVitie-Wilson
version [\mw] of Gale and Shapley's original stable marriage algorithm. The
basic idea is to look for cycles among the traders' best choices, and to put
such cycles into the allocation whenever they are found.

{\narrower\smallskip\noindent
{\bf A1.} [Initialize.] Set $(q_1,\ldots,q_n)\leftarrow (-1,\ldots,-1)$, 
$(r_1,\ldots,r_n)\leftarrow (0,\ldots,0)$, $(g_1,\ldots,g_n)\leftarrow
(0,\ldots,0)$, and $t\leftarrow 0$. (During this algorithm, $t$~will be a
trader who makes proposals to other traders, or zero when a new trader needs to
enter the picture. Variable~$q_k$ will be the number of a trader who currently
wants trader~$k$'s goods, or $q_k=0$ if trader~$k$ has expressed interest in
somebody else's wares but nobody has reciprocated; $q_k=-1$ if trader~$k$ has
not yet entered. Variable~~$r_k$ is the position of trader~$k$ in his list,
the number of proposals he has made.
Variable~$g_k$ is trader~$k$'s allocation, or 0 if no allocation has yet been
made.)

\smallskip\noindent
{\bf A2.} [Introduce a new trader.] (At this point $t=0$, and $g_k=0$ iff
$q_k=-1$. The traders with $g_k>0$ have been assigned a permutation of their
goods.) If all $g_k$ are nonzero, the algorithm terminates. Otherwise, set $t$
to some~$k$ with $g_k=0$, and set $q_t\leftarrow 0$.

\smallskip\noindent
{\bf A3.} [Propose.] Increase $r_t$ by 1, then set $s=p_{tr_t}$. If $g_s>0$,
repeat this step. (Trader~$t$ has expressed interest in his best remaining
choice,~$s$.) 

\smallskip\noindent
{\bf A4.} [Is $s$ spoken for?] If $q_s\geq 0$, go to step~A5. Otherwise set
$q_s\leftarrow t$ and $t\leftarrow s$, then return to~A3.

\smallskip\noindent
{\bf A5.} [Remove a cycle.] (There is now a cycle $s=s_1\rightarrow
s_2\rightarrow \cdots\rightarrow s_m=s$, where $s_{j+1}$ is the best remaining
choice of~$s_j$. This cycle must be part of any stable allocation, so we
incorporate it into~$g$.) Set $t\leftarrow q_s$; then repeatedly set
$g_s\leftarrow p_{sr_s}$ and $s\leftarrow p_{sr_s}$ until finding $g_s>0$. If
$t=0$, return to~A2, otherwise go to~A3.\quad\pfbox
\smallskip}

\smallskip
In step A3 there is always a path
$$0=t_0\rightarrow t_1\rightarrow t_2\rightarrow \cdots \rightarrow t_m=t$$
connecting all traders~$k$ such that $g_k=0$ and $q_k\geq 0$. Trader~$t_1$
entered in step~A2, and $t_{j+1}$ is the best remaining choice of~$t_j$, for
$1\leq j<m$; also $q_{t_j}=t_{j-1}$ for $1\leq j\leq m$. These invariant
relations justify the parenthesized assertions within the algorithm.

The final allocation $g_1\ldots g_n$ is stable. For if $h_k\leq g_k\ (k)$ for
all~$h_k$ in some coalition~$C$, and if some $h_k<g_k\ (k)$, then $h_k=p_{kr}$
for some $r<r_k$, so $h_k$ was rejected by the algorithm. 
When the algorithm changed~$r_k$ from~$r$ to $r+1$, it had already found $h_k$
to be the best remaining choice of some other trader~$s$, and it had assigned
$g_s=h_k$. Therefore $h_s>g_s\ (s)$. 

Moreover, the stable allocation is unique. If $g$ is assigned differently on
any cycle that leads to step~A5, that cycle will be a coalition violating~(2).

\meno
{\bf 3. A constructive lemma.}
We have observed that the stable allocation will be found by a
first-come-first-served algorithm equivalent to uniform hashing, using at least
one permutation~$\pi$ to give priority to the traders. 
For example, in the introduction we considered a case where $n=3$ and the
stable allocation was $1\,3\,2$. Any $\pi$ in which trader~3 has priority over
trader~1 will find this allocation.

We can also consider permutations of the preference lists. Let
$\sigma=\sigma(1)\ldots\sigma(n)$ be a permutation of $\{1,\ldots,n\}$, and
suppose that trader~$\sigma(k)$ uses list $p'_{\sigma(k)}=p_k$. This will
permute the locally optimum allocations, and it may also change the stable
allocation. For example, if $\sigma$ is $2\,1\,3$, so that
$$\eqalign{p'_2&=2\,1\,3\,,\cr
p'_1&=3\,1\,2\,,\cr
p'_3&=2\,3\,1\,,\cr}$$
the stable allocation becomes $g'_2g'_1g'_3=2\,1\,3$, because $g'_2$ must be~2
and then $g'_3$ must be~3. Allocating goods in the order $2\,1\,3$ was the
worst of the locally stable alternatives when $\sigma$ was the identity
permutation $1\,2\,3$, but it is best in the modified problem. Shuffling the
preference lists corresponds to shuffling the goods that the traders started
with.

We are now ready to prove a key fact about stable allocation. Let us say
that the prioritization~$\pi$ is {\it consistent with\/} the
shuffling~$\sigma$, 
with respect to preferences $p=p_1\ldots p_n$,
 if the locally optimum allocation
$g_1\ldots g_n$ obtained by uniform hashing with priorities~$\pi$ is the stable
allocation $g'_{\sigma(1)}\ldots g'_{\sigma(n)}$ when $p'_{\sigma(k)}=p_k$. For
example, $\pi(1)\,\pi(2)\,\pi(3)=1\,3\,2$ produces the locally optimum
$2\,1\,3$, so $\pi$ is consistent with the shuffling $\sigma(1)\,\sigma(2)\,
\sigma(3)=2\,1\,3$ just considered.

\proclaim
Lemma. Let $p$ be any sequence of preference lists. There is a one-to-one
correspondence between all permutations $\pi$ of $\{1,\ldots,n\}$ and all
permutations~$\sigma$ such that, if $\pi$ corresponds to~$\sigma$, the
prioritization~$\pi$ is consistent with the shuffling~$\sigma$.

\meno
{\bf Proof.}
Given $p$ and $\pi$, suppose uniform hashing with priorities~$\pi$ produces the
locally optimum allocation $g_1\ldots g_n$. Write the preference lists in rows,
with each~$g_k$ circled in its list~$p_k$. Delete all elements to the right
of~$g_k$.

We will construct a shuffling $\sigma$ whose stable allocation agrees with~$g$.
The construction involves two dynamically growing sets~$X$ and~$Y$, whose
significance will become clear momentarily. Initially 
$X\leftarrow\emptyset$ and $Y$
is the set of all~$k$ where $g_k=p_{k1}$ (i.e., all row numbers in which the
circled element is all by itself). Set $m\leftarrow 0$; as the construction
proceeds, we will have defined
$\sigma\bigl(\pi(1)\bigr),\ldots,\sigma\bigl(\pi(m)\bigr)$ as a permutation of
$\{g_{\pi(1)},\ldots,g_{\pi(m)}\}$, and we will have $X\subseteq Y$,
$\pi(m+1)\in Y\setminus X$.

Find the minimum $k$ in $m<k\leq n$ such that either $k=n$ or
$\bigl(\pi(k+1)\in Y\setminus X \hbox{ and } \pi(k+1)>\pi(m+1)\bigr)$ or
$\pi(k+1)\not\in Y$. Define $\sigma\bigl(\pi(j)\bigr)=g_{\pi(j+1)}$ for
$m<j<k$ and $\sigma\bigl(\pi(k)\bigr)=g_{\pi(m+1)}$. If $k=n$, the construction
is complete. Otherwise, remove $g_{\pi(m+1)},\ldots,g_{\pi(k)}$ from all
preference lists where they aren't circled. If $\pi(k+1)\not\in Y$, set
$X\leftarrow Y$ and let $Y$ be the set of all rows whose first elements are now
circled. (Since $\pi$ defines~$g$ by uniform hashing, $\pi(k+1)$ will be in the
new~$Y$.) Set $m\gets k$ and repeat the instructions of this paragraph.

A worked example will help clarify this construction.
Let
$$\pi(1)\ldots\pi(9)=5\,3\,4\,9\,1\,8\,2\,7\,6\,.$$
 Table 1 shows a sequence
of preference lists for $n=9$ in which $\pi$ defines the locally optimum
allocation indicated by circled elements. All elements to the right of the
circled ones have been erased, since they are irrelevant for our present
purposes.
$$\vcenter{\halign{#\hfil\quad&$#$\hfil\qquad
&#\hfil\quad&$#$\hfil\ &$#$\hfil\qquad
&#\hfil\quad&$#$\hfil\ &$#$\hfil\ &$#$\hfil\cr
1:&\cd 3&4:&1&\cd 5&7:&5&3&\cd 8\cr
\noalign{\smallskip}
2:&\cd4&5:&\cd9&&8:&9&\cd7\cr
\noalign{\smallskip} 
3:&\cd1&6:&\cd6&&9:&\cd2\cr}}$$
\centerline{Table 1. Preference lists and their stable allocation}

The stable allocation determined by these preference lists happens to coincide
with 
the circled elements in Table~1, so in this case 
the priorities~$\pi$ produce the stable
allocation; but our construction works for any~$\pi$, whether or not its
locally optimum allocation is stable. Initially $m=0$, $X=\emptyset$, and
$Y=\{1,2,3,5,6,9\}$. According to the rules stated, we proceed to
set $k=2$, since $\pi(3)\not\in Y$. So we define $\sigma(5)=g_3=1$ and
$\sigma(3)=9$; then we delete~1 and~9 from lists~4 and~8, and we set
$X\leftarrow \{1,2,3,5,6,9\}$, $Y\leftarrow\{1,2,3,4,5,6,8,9\}$, $m\leftarrow
2$. Next, $k=5$ since $\pi(6)=8\in Y\setminus X$ and $8>4=\pi(3)$. This time
$\sigma(4)=2$, $\sigma(9)=3$, $\sigma(1)=5$. We delete 2, 3, and~5 where they
are not circled. After setting $m\leftarrow 5$ we have $k=7$, because
$\pi(8)\not\in Y$. (Notice that $7\not\in Y$, even though row~7 now contains
only its circled element. The construction changes $Y$ only in the case
$\pi(k+1)\not\in Y$.) This time $\sigma(8)=4$, $\sigma(2)=7$, $X\leftarrow
\{1,2,3,4,5,6,8,9\}$, $Y\leftarrow\{1,\ldots,9\}$, $m\leftarrow 7$. On the
final round we set $\sigma(7)=6$ and $\sigma(6)=8$; the shuffled preference
lists are shown in Table~2. It is easy to verify that their stable allocation
matches that of Table~1, using the algorithm given earlier.
$$\vcenter{\halign{#\hfil\quad&$#$\hfil\qquad
&#\hfil\quad&$#$\hfil\ &$#$\hfil\qquad
&#\hfil\quad&$#$\hfil\ &$#$\hfil\ &$#$\hfil\cr
5:&\cd3&2:&1&\cd5&6:&5&3&\cd8\cr
\noalign{\smallskip}
7:&\cd4&1:&\cd9&&4:&9&\cd7\cr
\noalign{\smallskip} 
9:&\cd1&8:&\cd6&&3:&\cd2\cr}}$$
\centerline{Table 2. Shuffled precedence lists having the same stable
allocation as Table~1.}

The inverse construction is analogous. If $\sigma$ is any shuffling, circle its
stable allocation and prepare an array like Table~2. Begin with $m\leftarrow
0$, $X=\emptyset$, and $Y$ as before. Then repeatedly consider all cycles
$$\sigma(a_1)\leftarrow \sigma(a_2)\leftarrow
\cdots\leftarrow\sigma(a_t)\leftarrow\sigma(a_1)$$
formed by elements $\{a_1,\ldots,a_t\}\subseteq Y\setminus\{\pi(1),\ldots,
\pi(m)\}$. Here $\sigma(a)\leftarrow\sigma(b)$ means that $\sigma(a)$ is the
(circled) element in list $\sigma(b)$; for example,
$\sigma(4)\leftarrow\sigma(9)$ in Table~2 because the circled element in list
$\sigma(9)=3$ is $2=\sigma(4)$.
The properties of stable allocation guarantee that at least one such cycle
exists, and our inverse construction will guarantee that each cycle will
contain at least one $a_i\not\in X$. Call the largest such $a_i$ the {\it cycle
leader}, and renumber the subscripts so that $a_1$ is the cycle leader. Take
the cycle with smallest leader, and set $\pi(m+1)\leftarrow
a_1,\ldots,\pi(m+t)\leftarrow a_t$. 
Remove $\sigma(a_1),\ldots,\sigma(a_t)$ from the tableau in places where they
are not circled.
Then set $m\leftarrow m+t$ and repeat the same process until all cycles have
been recorded in~$\pi$. Then set $X\leftarrow Y$ and let $Y$ be the row numbers
that now have but a single element. Repetition of these steps will produce a
priority permutation~$\pi$ consistent with~$\sigma$. 

It is not difficult to verify that these constructions invert each other. The
reader will find easily, for example, that the permutation $\sigma(1)\ldots
\sigma(9)=5\,7\,9\,2\,1\,8\,6\,4\,3$ in Table~2 leads back to
$\pi(1)\ldots\pi(9)=5\,3\,4\,9\,1\,8\,2\,7\,6$. \ \pfbox

\meno
{\bf 4. A theorem.}
The lemma we have just proved makes it easy to establish the main result of
this note. We say that uniform hashing on $p_1\ldots p_n$ with priorities~$\pi$
produces ranks $r_1\ldots r_n$ if $r_{\pi(k)}$ is minimum such that
$$p_{\pi(k)r_{\pi(k)}}\not\in
\{p_{\pi(1)r_{\pi(1)}},\ldots,p_{\pi(k-1)r_{\pi(k-1)}}\}$$
for $1\leq k\leq n$.

\proclaim
Theorem. When preference lists $p_1\ldots p_n$ are independent and uniformly
random, the probability that the stable allocation $g_1\ldots
g_n=p_{1r_1}\ldots p_{nr_n}$ has a given value of the (unordered) multiset
$\{r_1,\ldots,r_n\}$ is the same as the probability that uniform hashing yields
$\{r_1,\ldots,r_n\}$.

\noindent
{\bf Proof.}
Let $\{r_1,\ldots,r_n\}$ be any given multiset. If $p=p_1\ldots p_n$ is any
sequence of preferences and $\sigma$ is any permutation of $\{1,\ldots,n\}$,
let $s(\sigma,p)=1$ if $\{r_1,\ldots,r_n\}$ is the multiset of ranks in the
stable allocation when trader $\sigma(k)$ has preference list~$p_k$; otherwise
$s(\sigma,p)=0$. 
Then the probability that stable allocation on random preferences has ranks
$\{r_1,\ldots,r_n\}$~is
$${1\over n!^n}\,\sum_p\,s(\sigma,p)$$
for any fixed $\sigma$. 

Similarly, if $\pi$ is any permutation of $\{1,\ldots,n\}$, let $h(\pi,p)=1$
iff $\{r_1,\ldots,r_n\}$ is the multiset of ranks produced by uniform hashing
with priorities~$\pi$.
Then the probability that uniform hashing on random preferences has ranks
$\{r_1,\ldots,r_n\}$~is
$${1\over n!^n}\,\sum_p\,h(\pi,p)$$
for any fixed $\pi$.

We want to show that these sums are equal. This is now obvious, because the
lemma implies that
$$\sum_{\sigma}\,\sum_ps(\sigma,p)=\sum_p\,\sum_{\sigma}s(\sigma,p)
=\sum_p\,\sum_{\pi}h(\pi,p)=\sum_{\pi}\,\sum_ph(\pi,p)$$
and we simply divide by $n!^{n+1}$. \ \pfbox

\medskip
Notice that this proof of the theorem remains valid even when the preference
lists $p_1\ldots p_n$ are not uniformly random. All we are assuming is a
symmetry condition, that shuffled preference lists $p_{\sigma(1)}\ldots
p_{\sigma(n)}$ have the same distribution for all~$\sigma$.

\meno
{\bf 5. Corollaries.}
The analysis of uniform hashing is quite simple, so our theorem immediately
characterizes many properties of the ranks in random stable allocations. For
example, let us find the expected value of
$$(z+r_1)(z+r_2)\ldots(z+r_n)\,;$$
this polynomial is clearly a function of the multiset $\{r_1,\ldots,r_n\}$, so
we can analyze it by considering its behavior with respect to uniform hashing.

Let $q_{kj}$ be the probability that $r_k>j$ in uniform hashing. This is the
probability that
$p_{k1},\ldots, p_{kj}\in\{p_{1r_1},\ldots,p_{(k-1)r_{k-1}}\}$,~so
$$q_{kj}=\left({k-1\over n}\right)\left({k-2\over n-1}\right)\ldots
\left({k-j\over n-j+1}\right)={k-1\choose j}
  \left/{n\choose j}\right.\,.\eqno(1)$$
Standard binomial coefficient summation techniques [\cm] show that
$$\sum_{j=0}^{\infty}{j\choose m}q_{kj}={n+1\over n+m+2-k}\,\left.
{k-1\choose m}\right/{n+m+1-k\choose m}\,.\eqno(2)$$
The expected value of $(z+r_1)\ldots (z+r_n)$ is therefore
$$\eqalignno{%
\sum_{r_1,\ldots,r_n}&\,\prod_{k=1}^n\bigl(q_{k(r_k-1)}-q_{kr_k}\bigr)
(z+r_k)\cr
\noalign{\smallskip}
&=\prod_{k=1}^n\,\sum_{r=1}^{\infty}\bigl(q_{k(r-1)}-q_{kr}\bigr)(z+r)\cr
\noalign{\smallskip}
&=\prod_{k=1}^m\left(z+\sum_{j=0}^{\infty}q_{kj}\right)\cr
\noalign{\smallskip}
&=\prod_{k=1}^n\left(z+{n+1\over n+2-k}\right)={1\over
(n+1)!}\,\prod_{k=2}^{n+1} (kz+n+1)\,.&(3)\cr}$$
In particular, the expected value of $r_1+\cdots +r_n$, which is the
coefficient of~$z^{n-1}$,~is
$$\sum_{k=1}^n\,{n+1\over n+2-k}=(n+1)(H_{n+1}-1)=(n+1)H_n-n\,.\eqno(4)$$
The other coefficients can be expressed in terms of Stirling cycle numbers if
we note that
$$E(z+r_1)\ldots (z+r_n)(z+n+1)={1\over (n+1)!}\,\prod_{k=1}^{n+1}(kz+n+1)
=\sum_{k=0}^{n+1}\,{\,n+2\,\brack k+1}\,{(n+1)^k\over
(n+1)!}\,z^{n+1-k}\,.\eqno(5)$$
For example, the coefficient of $z^{n-2}$ in $E(z+r_1)\ldots (z+r_n)$~is
$$\eqalignno{{\,n+2\,\brack 3}\,&{(n+1)^2\over (n+1)!}-(n+1){\,n+2\,\brack 2}\,
{(n+1)\over (n+1)!}+(n+1)^2{\,n+2\,\brack 1}\,{1\over (n+1)!}\cr
\noalign{\smallskip}
&=(n+1)^2\left({\textstyle{1\over 2}}\bigl(H^2_{n+1}-H^{(2)}_{n+1}\bigr)
-H_{n+1}+1\right)\cr
\noalign{\smallskip}
&={(n+1)^2\over
2}\bigl(H_n^2-H_n^{(2)}\bigr)-n(n+1)(H_n-1)\,;&(6)\cr}$$
see [\cm, exercise 6.33].

So far we have used only the case $m=0$ of (2). A similar argument, using
$m=1$, shows that
$$E(z+r_1^2)\ldots (z+r_n^2)=\prod_{k=1}^n\,\left(z+{(n+1)(n+1+k)\over
(n+2-k)(n+3-k)}\right)\,.\eqno(7)$$
In particular,
$$\eqalignno{E(r_1^2+\cdots r_n^2)&=\sum_{k=1}^n\,
{(n+1)(n+1+k)\over (n+2-k)(n+3-k)}\cr
\noalign{\smallskip}
&=(n+1)\,\sum_{k=1}^n\,
\left({2n+4\over (n+2-k)(n+3-k)}-{1\over n+2-k}\right)\cr
\noalign{\smallskip}
&=(n+1)(n-H_{n+1}+1)=(n+1)(n-H_n)+n\,.&(8)\cr}$$
Hence, by (6) and (8),
$$\eqalignno{E\bigl((r_1+\cdots +r_n)^2\bigr)
&=E(r_1^2+\cdots +r_n^2)+2\,[z^{n-2}]\,E(z+r_1)\ldots (z+r_n)\cr
\noalign{\smallskip}
&=(n+1)^2\bigl(H_n^2-H_n^{(2)}\bigr)-(n+1)(2n+1)H_n+n(3n+4)\,.&(9)\cr}$$
The expected value of the variance of the ranks is therefore
$$E\left({r_1^2+\cdots +r_n^2\over n}\right)-E\biggl(\left({r_1+\cdots
+r_n\over n}\right)^{\!2}\biggr)=n+O(\log n)^2\eqno(10)$$
while the variance of the rank sum is
$$\eqalignno{E\bigl((r_1\cdots +r_n)^2\bigr)-\bigl(E(r_1+\cdots +r_n)\bigr)^2
&=2n(n+2)-(n+1)^2H_n^{(2)}-(n+1)H_n\cr
\noalign{\smallskip}
&=\left(2-{\pi^2\over 6}\right)\,n^2+O(n\log n)\,.&(11)\cr}$$

The final rank $r_n$ in uniform hashing is uniformly distributed in
$\{1,\ldots,n\}$. Therefore the probability is $\geq {1\over 2}$ that at least
one trader in a random stable allocation will have rank $\geq {1\over 2}n$.
(He will be left with 
a~piece of goods he doesn't like very much, since it's in the lower half
of his list.) Indeed, the probability that
$\max(r_1,\ldots,r_n)\leq {1\over 2}n$ is exactly
$$(1-q_{1m})(1-q_{2m})\ldots(1-q_{nm})\,,$$
where $m=\lfloor{1\over 2}n\rfloor$; this is asymptotically
$${\textstyle{(1-{1\over 2})(1-{1\over 4})(1-{1\over 8})(1-{1\over 16})}}
\ldots \approx .288788\,.\eqno(12)$$

\meno
{\bf 6. Conclusions and conjectures.}
The running time of the
 simple algorithm we have presented for stable allocation is essentially
proportional to the sum of ranks in the unique allocation, 
$r_1+\cdots +r_n$. We have proved that the statistical
properties of any symmetric function
of $(r_1\ldots r_n)$ are identical to the corresponding statistics for uniform
hashing, provided only that the distribution of preference lists $p_1\ldots p_n$
is invariant under shuffling. When the preferences are uniformly random, the
expected value of $r_1+\cdots +r_n$ is exactly $(n+1)H_n-n$, and the standard
deviation is $O(n)$.

Uniform hashing is equivalent to the classical stable marriage problem when all
the girls have the same preference list. (See [\ms, pages 65--67].) 
Perhaps it is worthwhile to repeat here the main research problem about stable
marriages that was advertised in the author's lectures of 1975 [\ms] and not
yet resolved: If the girls have any fixed set of preferences and the boys
propose at random, is the expected rank sum $r_1+\cdots +r_n$ of the
male-optimum stable marriage always $\geq (n+1)H_n-n$? In other words, does the
case of equal preferences for the girls (uniform hashing) give the greatest
lower bound for $E(r_1+\cdots +r_0)$? If so, the average would be tightly
bounded, because the upper bound $(n-1)H_n+1$ is easy to prove [\ms, page~43]. 

In fact, computer experiments for small~$n$ suggest that the 
maximum value of $E(r_1+\cdots
+r_n)$, when the girls have a fixed set of preferences and the boys propose
independently at random, is obtained if and only if the girls' preferences are
cyclic, in the sense that we could rename boys and girls so that girl~$j$'s
$k$\/th~choice is congruent to $j+k$ (mod $n$).

Both conjectures about min and max $E(r_1+\cdots +r_n)$ have been verified by
exhaustive enumeration when $n\leq 4$, and in several hundred
random experiments when $n=5$.
Presumably there is a (simple?) way to prove that, in some sense, the more the
girls agree in their ranking, the less the men will have to propose, on
the average.

Is there a simple expression for $E(r_1+\cdots +r_n)$ when the girls'
preferences are cyclic? For $n=3,4,5$ the values are respectively $306/3!^3$,
$884224/4!^4$, $104035560000/5!^5$. When $n=4$, the worst seven preference
matrices for the girls are
$$\tabskip\centering
\halign to\displaywidth{\hfil$#$\tabskip2em&\hfil$#$&\hfil$#$&\hfil$#$%
&\hfil$#$&\hfil$#$&\hfil$#$\tabskip\centering\cr
1\,2\,3\,4&1\,2\,3\,4&1\,2\,3\,4&1\,2\,3\,4&1\,2\,3\,4&1\,2\,3\,4&1\,2\,3\,4\cr
2\,3\,4\,1&2\,3\,1\,4&2\,3\,1\,4&2\,3\,4\,1&2\,3\,1\,4&2\,3\,1\,4&1\,2\,3\,4\cr
3\,4\,1\,2&3\,4\,1\,2&3\,4\,2\,1&3\,1\,4\,2&3\,4\,1\,2&3\,4\,1\,2&3\,4\,1\,2\cr
4\,1\,2\,3&4\,1\,2\,3&4\,1\,2\,3&4\,1\,2\,3&4\,1\,3\,2&4\,2\,1\,3&4\,2\,3\,1\cr
\noalign{\medskip\hbox{with respective total rank sums}\medskip}
884224\rlap,&879488\rlap,&875264\rlap,&875072\rlap,&874752\rlap,%
&874624\rlap,&872192\rlap.\cr}$$
All preference matrices not isomorphic to these seven,
 under renumbering of boys and girls, have smaller total rank sum over all
$4!^4$ preference matrices for the boys.

\meno
{\bf Acknowledgment.}
I want to thank Boris Pittel for introducing me to this problem and for
patiently correcting my original misunderstanding of the definitions.

\medskip
\centerline{\bf References}
\medskip

\bib 
[\fp]
Alan M. Frieze and Boris G. Pittel, ``Probabilistic analysis of an algorithm in
the theory of markets in indivisible goods,'' {\sl Annals of Applied
Probability\/ \bf5} (1995), 768--808.

\smallskip
\bib
[\gs]
D. Gale and L. S. Shapley, ``College admissions and the stability of
marriage,'' {\sl American Mathematical Monthly\/ \bf 69} (1962), 9--15.

\smallskip
\bib
[\cm]
Ronald L. Graham, Donald E. Knuth, and Oren Patashnik, {\sl Concrete
Mathematics\/} (Reading, Massachusetts: \AW, 1989).

\smallskip
\bib
[\gi]
Dan Gusfield and Robert W. Irving, {\sl The Stable Marriage Problem\/}
(Cambridge, Mass.: MIT Press, 1989).

\smallskip
\bib
[\ks]
Donald E. Knuth, {\sl Sorting and Searching\/} (Reading, Massachusetts:
\AW, 1973).

\smallskip
\bib
[\ms]
Donald E. Knuth, {\sl Mariages Stables\/} (Montr\'eal: Les Presses de
l'Universit\'e de Montr\'eal, 1976).

\smallskip
\bib
[\mw]
D. G. McVitie and L. B. Wilson, ``The stable marriage problem,'' {\sl
Communications of the ACM\/ \bf 14} (1971), 486--492.

\smallskip
\bib
[\ss]
L. S. Shapley and H. Scarf, ``On cores and indivisibility,'' {\sl Journal of
Mathematical Economics\/ \bf 1} (1974), 23--38.

\bye